\documentclass[12pt]{article}
\usepackage{amsmath,amssymb,amsthm,amsfonts,hhline,color}
\usepackage{bbm}
\usepackage[titletoc,title]{appendix}
\usepackage[hypertex]{hyperref}
\usepackage{comment}
\usepackage{color}
\usepackage{graphicx}
\usepackage{mathrsfs}

\newcommand{\ds}{\displaystyle}

\newcommand{\dfr}[2]{\dfrac{#1}{#2}}
\newcommand{\cd}{\cdot}
\newcommand{\cds}{\cdots}
\newcommand{\dsum}{\displaystyle \sum}

\renewcommand{\l}{\left}
\renewcommand{\r}{\right}

\newcommand{\q}{\quad}
\newcommand{\qq}{\qquad}

\newcommand{\la}{\langle}
\newcommand{\ra}{\rangle}

\newcommand{\abs}[1]{\lvert{#1}\rvert}

\newcommand{\Z}{\mathbb{Z}}
\newcommand{\C}{\mathbb{C}}
\newcommand{\R}{\mathbb{R}}
\newcommand{\N}{\mathbb{N}}
\newcommand{\Q}{\mathbb{Q}}

\newcommand{\vir}{\mathrm{Vir}}
\newcommand{\aut}{\mathrm{Aut}}
\newcommand{\wt}{\mathrm{wt}}
\newcommand{\tr}{\mathrm{Tr}}

\newcommand{\id}{\mathrm{id}}

\newcommand{\ad}{\mathrm{ad}}

\newcommand{\Span}{\mathrm{Span}}

\newcommand{\w}{\omega}
\newcommand{\vac}{\mathbbm{1}}

\newcommand{\mymid}{\,|\,}

\makeatletter
\@addtoreset{equation}{section}

\makeatother

\theoremstyle{plain}

\newtheorem{thm}{Theorem}[section]
\newtheorem{prop}[thm]{Proposition}
\newtheorem{lem}[thm]{Lemma}
\newtheorem{cor}[thm]{Corollary}

\theoremstyle{definition}
\newtheorem{df}[thm]{Definition}
\newtheorem{conj}[thm]{Conjecture}

\newtheorem{cond}{Condition}

\newtheorem{rem}[thm]{Remark}

\newcommand{\sfr}[2]{\leavevmode\kern-.05em
  \raise.5ex\hbox{\the\scriptfont0 #1}\kern-.1em
  /\kern-.15em\lower.25ex\hbox{\the\scriptfont0 #2}\kern.02em}

\newcommand{\sshf}{\leavevmode\kern-.05em
  \raise.3ex\hbox{\the\scriptfont0 1}\kern-.1em
  /\kern-.10em\lower.5ex\hbox{\the\scriptfont0 2}\kern.02em}

\newcommand{\shf}{\sfr{1}{2}}
\DeclareMathOperator*{\tensor}{\otimes}
\DeclareMathOperator*{\fusion}{\boxtimes}

\newcommand{\pf}{\noindent {\bf Proof:}\q }

\newcommand{\com}{\mathrm{Com}}
\newcommand{\sym}{\mathfrak{S}}

\renewcommand{\ker}{\mathrm{Ker}}
\renewcommand{\o}{\mathrm{o}}

\newcommand{\A}{\mathscr{A}}

\title{Vertex operator algebras generated by \\ Ising vectors of $\sigma$-type}

\author{
  Cuipo  Jiang\thanks{Partially supported by China NSF grants 11771281, 11531004.}
  \medskip\\
  \textit{\small School of Mathematical Science,  Shanghai Jiao Tong University} \\
  \textit{\small Shanghai, 200240, China}
  \\
  {\small email: \texttt{cpjiang@sjtu.edu.cn}}
  \bigskip\\
  Ching Hung Lam\thanks{Partially supported by MoST grant 104-2115-M-001-004-MY3 of Taiwan.}%
  \medskip\\
  \textit{\small Institute of Mathematics, Academia Sinica, Taipei, Taiwan 10617}\\
  {\small e-mail: \texttt{chlam@math.sinica.edu.tw}}
  \bigskip\\
  Hiroshi Yamauchi\footnote{Partially supported by JSPS KAKENHI Grant Number
  JP16K05073. }
  \medskip\\
  \textit{\small Department of Mathematics,  Tokyo Woman's Christian University}\\
  \textit{\small 2-6-1 Zempukuji, Suginami-ku, Tokyo 167-8585, Japan}\\
  {\small e-mail: \texttt{yamauchi@lab.twcu.ac.jp}}
}
\date{}

\begin{document}


\maketitle

\begin{abstract}
We prove the uniqueness of the simple vertex operator algebra of OZ-type
generated by Ising vectors of $\sigma$-type.
We also prove that the simplicity can be omitted if the Griess algebra is isomorphic to
the Matsuo algebra associated with the root system of type $A_n$.
\end{abstract}

\pagestyle{plain}

\baselineskip 6mm


\section{Introduction}

An \emph{Ising vector} is a Virasoro vector of a vertex operator algebra (VOA)
generating a simple Virasoro subVOA of central charge 1/2.
Vertex operator algebras generated by Ising vectors are interesting from the point of view of 
(finite) group theory since each Ising vector defines an involution called Miyamoto involution
acting on them (cf.~\cite{Mi}).
The group generated by Miyamoto involutions associated with Ising vectors of $\sigma$-type
(see Section 3 for definition) forms a 3-transposition group.
When the VOAs have compact real forms, the 3-transposition groups obtained in this manner are
of symplectic type (cf.~\cite{CH}) and are completely classified in \cite{Ma}.
In \cite{Ma}, a complete list of 3-transposition groups generated by Miyamoto involutions
associated with Ising vectors of $\sigma$-type as well as examples of VOAs realizing those
groups are presented.
The purpose of the paper \cite{Ma} is mainly on the classification of 3-transposition groups
but not of VOAs.
The classification of VOAs realizing 3-transposition groups has not been established so far.
In this paper, we will prove the uniqueness of simple VOA structure of a VOA generated by
Ising vectors of $\sigma$-type when it is of OZ-type, i.e., it is of CFT-type and the
weight one subspace is trivial.
Our result complements Matsuo's classification in the sense that the examples of VOAs
in \cite{Ma} are indeed unique under the conditions considered there.

Let us explain our results precisely.
Let $V$ be a VOA of OZ-type and $E_V$ the set of Ising vectors of $V$ of $\sigma$-type.
For $e\in E_V$, one can define the Miyamoto involution $\sigma_e$ acting on $V$.
It is shown in \cite{Mi} that the group $G_V$ generated by $\{ \sigma_e \mymid e\in E_V\}$
is a 3-transposition group.
The linear span of $E_V$ forms a subalgebra of the Griess algebra of $V$ which has a
description in terms of the 3-transposition group $G_V$.
This algebra structure is called the Matsuo algebra $B_{\sshf,\sshf}(G_V)$ associated with $G_V$
(cf.~\cite{HRS,Ma}).
Suppose $V$ is generated by $E_V$ as a VOA.
The Griess algebra of $V$ coincides with a homomorphic image of the Matsuo algebra
$B_{\sshf,\sshf}(G_V)$ and $G_V$ is a center-free 3-transposition group in $\aut(V)$ 
(cf.~\cite{Mi,Ma}).
We will prove that the bilinear form on $V$ is uniquely determined by its Griess algebra
structure in Proposition \ref{prop:3.7}.
Since $V$ is assumed to be of OZ-type, it has a unique simple quotient given by
the non-degenerate quotient with respect to the bilinear form on $V$.
Since the radical of the bilinear form on $V$ is uniquely determined by the Griess algebra,
we will prove in Theorem \ref{thm:3.6} the uniqueness of the VOA structure of $V$ if it is simple.
When $V$ is simple, the Griess algebra of $V$ is isomorphic to the quotient of Matsuo algebra
by the radical of the bilinear form and uniquely determined by $G_V$.
So our result shows that $V$ is uniquely determined by the 3-transposition group $G_V$
if $V$ is simple and $G_V$ is realized by Miyamoto involutions associated with Ising vectors
of $\sigma$-type.
In the case that $G_V$ is isomorphic to the symmetric group $\sym_{n}$ of degree $n$,
the Matsuo algebra $B_{\sshf,\sshf}(\sym_{n})$ is already non-degenerate and we will prove
in Theorem \ref{thm:4.6} that $V$ is uniquely determined without assuming its simplicity.
The simple VOA with $G_V=\sym_{n}$ is well studied in \cite{LS,JL} and our result gives a
characterization of this VOA in terms of the Griess algebra.

As a by-product, we can weaken the assumption on the positivity in \cite{Ma}.
By the uniqueness for the VOA corresponding to $\sym_3$,
we can prove in Proposition \ref{prop:3.13} that the 3-transposition group $G_V$ generated
by Ising vectors of $\sigma$-type is always of symplectic type without assuming that
$V$ has a compact real form.
However, since the existence of Ising vectors already requires the positivity on local
subalgebras, it seems that all VOAs satisfying our conditions have compact real forms and
we conjecture that Matsuo's list in \cite{Ma} is complete without assuming the positivity
(cf.~Conjecture \ref{conj:3.16}).

The organization of this article is as follows.
In Section 2, we review Matsuo algebras associated with 3-transposition groups
in a general setting.
In Section 3,  we will prove the uniqueness of simple VOA structure of a VOA under the
assumptions that it is of OZ-type and generated by Ising vectors of $\sigma$-type.
In Section 4,  we consider a VOA of OZ-type which is generated by Ising vectors of $\sigma$-type
and whose Griess algebra is isomorphic to the Matsuo algebra associated with the symmetric group
$\sym_{n+1}$ (or the Weyl group of type $A_n$).
We will prove the uniqueness of the VOA structure of such a VOA without assuming the simplicity.
The whole arguments in Section 4 are based on induction on $n$.

\paragraph{Acknowledgement.}
Part of this work has been done while the second and the third named authors were staying
at Shanghai Jiao Tong University in June 2016 and in November 2017.
They gratefully acknowledge the hospitality there.
The authors thank the referee for careful reading of the manuscript and valuable suggestions 
which greatly improved this article.
H.Y. thanks Atsushi Matsuo for valuable comments and helpful discussions.

\paragraph{Notation and terminology.}
In this paper, vertex operator algebras (VOAs) are defined over the complex number field $\C$
unless otherwise stated.
A VOA $V$ is \emph{of CFT-type} if it has the $L(0)$-grading $V=\oplus_{n\geq 0}V_n$
such that $V_0=\C \vac$, and is \emph{of One-Zero type} or simply \emph{of OZ-type}
if it is of CFT-type and $V_1=0$.
In this case, $V$ is equipped with a unique invariant bilinear form such that
$(\vac|\vac)=1$ (cf.~\cite{Li}).
In this paper, we only consider VOAs of OZ-type.
A real form $V_\R$ of $V$ is called \emph{compact} if the associated bilinear form
is positive definite.
For a subset $A$ of $V$, the subalgebra generated by $A$ is denoted by $\la A\ra$.
For $a\in V_n$, we define its weight by $\wt(a)=n$.
We write $Y(a,z)=\sum_{n\in \Z}a_{(n)}z^{-n-1}$ for $a\in V$ and define
its \emph{zero-mode} by $\o(a):=a_{(\wt(a)-1)}$ if $a$ is homogeneous, and extend it linearly.
The weight two subspace $V_2$ carries a structure of a commutative algebra defined
by the product $\o(a)b=a_{(1)}b$ for $a$, $b\in V_2$.
This algebra is called the \emph{Griess algebra} of $V$.
A \emph{Virasoro vector} is an element $e\in V_2$ such that $e_{(1)}e=2e$.
It is known (cf.~\cite{Mi}) that $L^e(n):=e_{(n+1)}$, $n\in\Z$ generate a Virasoro algebra of central
charge $c_e=2(e\mymid e)$ and hence the subalgebra $\la e\ra$ is isomorphic to a Virasoro VOA.
We denote by $\vir(e)$ the Virasoro algebra generated by $e$.
$\vir(e)$ has a standard triangular decomposition by 
$\vir(e)_\pm =\oplus_{n>0} \C L^e(\pm n)$ and $\vir(e)_0=\C L^e(0)\oplus \C c_e$.
The universal enveloping algebra of $\vir(e)$ is denoted by $\mathrm{U}(\vir(e))$
which is naturally graded by $\Z$ so that 
$\mathrm{U}(\vir(e))[j]=\{ x\in \mathrm{U}(\vir(e)) \mid [L^e(0),x]=j x\}$ for $j\in\Z$.
Two Virasoro vectors $a$ and $b$ are said to be \emph{orthogonal} if $a_{(n)}b=0$ for $n\geq 0$.
We denote by $L(c,h)$ the irreducible highest weight module over the Virasoro algebra
with the central charge $c$ and the highest weight $h$.
A \emph{simple} $c=c_e$ Virasoro vector $e\in V$ is a Virasoro vector such that
$\la e\ra \cong L(c_e,0)$.
A Virasoro vector $\w$ is called the \emph{conformal vector} of $V$ if each graded
subspace $V_n$ agrees with $\ker_V\, (\o(\w)-n)$ and satisfies
$\w_{(0)}a =a_{(-2)}\vac$ for all $a\in V$.
The half of the conformal vector gives the unit of the Griess algebra and
hence is uniquely determined.
We write $\w_{(n+1)}=L(n)$ for $n\in \Z$.
A \emph{subVOA} $(W,e)$ of $V$ is a pair of a subalgebra $W$ of $V$
together with a Virasoro vector $e$ in $W$ such that $e$ is the conformal vector of $W$.
We often omit $e$ and simply denote it by $W$.
A subVOA $W$ of $V$ is called \emph{full\,} or \emph{conformal} if $V$ and $W$ share 
the same conformal vector.
The commutant subalgebra of $(W,e)$ in $V$ is defined by $\com_V W:=\ker_V\,e_{(0)}$
(cf.~\cite{FZ}).

\section{Matsuo algebra of a 3-transposition group}

We recall the definition of 3-transposition groups.

\begin{df}
  A \emph{$3$-transposition group} is a pair $(G,I)$ of a group $G$ and a set of involutions
  $I$ of $G$ satisfying the following conditions.
  \\
  (1)~$G$ is generated by $I$.
  \\
  (2)~$I$ is closed under conjugations, i.e., if $a$, $b \in I$ then $a^b=bab\in I$.
  \\
  (3)~For any $a$ and $b\in I$, the order of $ab$ is bounded by 3.
  \\
  A 3-transposition group $(G,I)$ is called \emph{indecomposable} if $I$ is a conjugacy class
  of $G$.
  An indecomposable $(G,I)$ is called \emph{non-trivial} if $I$ is not a singleton, i.e.,
  $G$ is not cyclic.
\end{df}

Let $(G,I)$ be a 3-transposition group.
We define a graph structure on $I$ as follows: for $a$, $b\in I$, $a\sim b$ if and only if $a$ and $b$ are non-commutative.
If $a\sim b$ then the order of $ab$ is three and we have $a^b=b^a\in I$.
We set $a\circ b:=a^b=b^a$ if $a\sim b$.
It is clear that $I$ is a connected graph if and only if $I$ is a single conjugacy class of $G$.

Let $\alpha$, $\beta$ be non-zero complex numbers.
Let $B_{\alpha,\beta}(G,I)=\oplus_{i\in I} \C x^i$ be the vector space spanned by a formal basis
$\{ x^i\mid i\in I\}$ indexed by the set of involutions.
We define a bilinear product and a bilinear form on $B_{\alpha,\beta}(G,I)$ by
\begin{equation}\label{eq:2.1}
  x^ix^j:=
  \begin{cases}
    2x^i & \mbox{if $i=j$},
    \\
    \dfr{\alpha}{2} (x^i+x^j-x^{i\circ j}) & \mbox{if $i\sim j$},
    \\
    0 & \mbox{otherwise},
  \end{cases}
  ~~~~~~
  (x^i|x^j):=
  \begin{cases}
    ~\dfr{\beta}{2} & \mbox{if $i=j$},
    \medskip\\
    \dfr{\alpha\beta}{8} & \mbox{$i \sim j$},
    \medskip\\
    ~\,0 & \mbox{otherwise}.
  \end{cases}
\end{equation}
Then $B_{\alpha,\beta}(G,I)$ is a commutative non-associative algebra with
a symmetric invariant bilinear form (cf.~\cite{Ma}).
This algebra is called the \emph{Matsuo algebra} associated with a 3-transposition group $(G,I)$
with the accessory parameters $\alpha$ and $\beta$ (cf.~\cite{HRS}\footnote{%
Strictly speaking, a Matsuo algebra is associated with a partial triple system and
we are considering Matsuo algebras associated with the Fischer spaces of
3-transposition groups.}).
In the following, we will simply denote $B_{\alpha,\beta}(G,I)$ by $B_{\alpha,\beta}(G)$.

\begin{rem}
  For a subfield $\Q(\alpha,\beta)\subset K \subset \C$, we can define a canonical
  $K$-form $B_{\alpha,\beta}(G)_K=\oplus_{i\in I}Kx^i$ of $B_{\alpha,\beta}(G)$.
  In particular, if $\alpha$, $\beta$ are real, we can define a canonical
  $\R$-form $B_{\alpha,\beta}(G)_\R=\oplus_{i\in I}\R x^i$.
\end{rem}

Suppose $I$ is not a single conjugacy class.
Then there is a partition $I=I_1\sqcup I_2$ such that both $I_1$ and $I_2$ are closed under
conjugations.
Then $G_1=\la I_1\ra$ and $G_2=\la I_2\ra$ are 3-transposition subgroups of $G$.
Since $[G_1,G_2]=1$ and $G=G_1G_2$, $G$ is a central product of $G_1$ and $G_2$.
In this case, $B_{\alpha,\beta}(G)$ is a direct sum $B_{\alpha,\beta}(G_1)\oplus B_{\alpha,\beta}(G_2)$
of two-sided ideals which is also an orthogonal sum with respect to the bilinear form.
Therefore, the structure of $B_{\alpha,\beta}(G)$ is determined by its indecomposable components.

\begin{df}
  The radical of the bilinear form on $B_{\alpha,\beta}(G)$ forms an ideal.
  We call the quotient algebra of $B_{\alpha,\beta}(G)$ by the radical of
  the bilinear form the \emph{non-degenerate quotient}.
\end{df}

It is obvious that $G$ acts on $B_{\alpha,\beta}(G)$ by conjugation.
Namely, we can define a group homomorphism
\begin{equation}\label{eq:2.2}
  \rho: G \to \aut\, B_{\alpha,\beta}(G)
\end{equation}
by, for $i\in I$, letting $\rho_ix^j = x^{i\circ j}$ if $i \sim j$ and $\rho_i x^j=x^j$ otherwise.
It is shown in \cite{Ma} that $\rho$ is injective if and only if $G$ is non-trivial and
center-free.

Suppose $\alpha \ne 2$.
Then the adjoint action of $x^i$ is semisimple with eigenvalues contained in $\{2,0,\alpha\}$.

For $i\in I$, let
\[
  B_{\alpha,\beta}(G)[i;h]=\ker_{B_{\alpha,\beta}(G)}(\ad(x^i)-h)~~~\mbox{and}~~~
  B_{\alpha,\beta}^{i\pm}(G)=\{ v\in B_{\alpha,\beta}(G) \mid \rho_iv=\pm v\}.
\]
Then $B_{\alpha,\beta}(G)[i;2]= \C x^i$ and we have
\[
  B_{\alpha,\beta}^{i+}(G)= \C x^i \oplus B_{\alpha,\beta}(G)[i;0], ~~~
  B_{\alpha,\beta}^{i-}(G)= B_{\alpha,\beta}(G)[i;\alpha].
\]
Namely, the action of $i\in I$ on $B_{\alpha,\beta}(G)$ can be described by the adjoint action of
$x^i$.

Suppose $G$ is indecomposable.
Then the number $\abs{\{ j\in I \mid j\sim i\}}$ is independent of the choice of $i\in I$
if $I$ is finite.
We denote this number by $k$.
One can verify that
\[
  \l(\dsum_{i\in I} x^i\r) \cd x^j=\l(\dfr{k\alpha}2+2\r) x^j.
\]
If $k\alpha+4$ is non-zero,  then
\begin{equation}\label{eq:2.3}
  \w:=\dfr{4}{k\alpha+4}\dsum_{i\in I} x^i
\end{equation}
satisfies $\w v=2v$ for $v\in B_{\alpha,\beta}(G)$ and $\w$ gives twice the unity of
$B_{\alpha,\beta}(G)$.
By the invariance, one has $(x^i|\w)=(x^i|x^i)$ and $(\w|\w)=\dfr{2\beta\abs{I}}{k\alpha+4}$.

\begin{rem}
  In VOA theory, a Matsuo algebra corresponds to a Griess algebra generated by
  $c=\beta$ Virasoro vectors having two highest weights $0$ and $\alpha$ with
  binary fusion rules (cf.~\cite{Ma}).
  The vector $\w$ is the conformal vector of such a VOA and $2(\w|\w)$ gives the central charge.
\end{rem}

\section{VOAs generated by Ising vectors of $\sigma$-type}

Recall the unitary series of the Virasoro VOAs.
Let
\begin{equation}\label{eq:3.1}
\begin{array}{rl}
  c_n &:= 1-\dfr{6}{(n+2)(n+3)},\qq n=1,2,3,\dots ,
  \medskip\\
  h_{r,s}^{(n)} &:= \dfr{\l(r(n+3)-s(n+2)\r)^2-1}{4(n+2)(n+3)},~~
  1\leq r\leq n+1,~~1\leq s\leq n+2.
\end{array}
\end{equation}
It is shown in \cite{W} that $L(c_n,0)$ is rational and $L(c_n,h_{r,s}^{(n)})$,
$1\leq s\leq r\leq n+1$, are all the irreducible $L(c_n,0)$-modules
(see also \cite{DMZ}).
The fusion rules among $L(c_n,0)$-modules are computed in \cite{W} and given by
\begin{equation}\label{eq:3.2}
  L\l( c_n,h^{(n)}_{r_1,s_1}\r)\fusion L\l( c_n,h^{(n)}_{r_2,s_2}\r)
  = \dsum_{1\leq i\leq M\atop 1\leq j\leq N}
    L(c_n,h^{(n)}_{\abs{r_1-r_2}+2i-1,\abs{s_1-s_2}+2j-1}),
\end{equation}
where $M=\min \{ r_1,\,r_2,\,n+2-r_1,\,n+2-r_2\}$ and $N=\min \{ s_1,\,s_2,\,n+3-s_1,\,n+3-s_2\}$.

\begin{df}
  A Virasoro vector $e\in V_2$ with central charge $c$ is called {\it simple} if 
  $\la e\ra\cong L(c,0)$.
  A simple $c=1/2$ Virasoro vector is called an {\it Ising vector}.
\end{df}

It is well-known (cf.~\cite{DMZ}) that $L(\shf,0)$ is rational and has three irreducible modules,
$L(\shf,0)$, $L(\shf,\shf)$ and $L(\shf,\sfr{1}{16})$.


\begin{lem}\label{lem:3.8}
  Let $e\in V$ be an Ising vector and $v\in V$ a highest weight vector for $\vir(e)$ 
  with the highest weight $h\in \{ 0,1/2\}$.
  Then the following relations hold.
  \[
    {e_{(0)}} v =0 ~~\mbox{if~ $h=0$},~~~~~
    {e_{(0)}}^2\, v 
    =\dfr{4}{3}e_{(-1)}v ~~\mbox{if~ $h=1/2$}.
  \]
\end{lem}

\pf
Since $e_{(0)}v$ (resp.~$3{e_{(0)}}^2v - 4 e_{(-1)}v$) is a singular vector for $h=0$ 
(resp.~$h=1/2$), the lemma follows (cf.~\cite{KR}). 
\qed

The following lemma can be proved by using the fermionic construction of the SVOA $L(\shf,0)\oplus L(\shf,\shf)$. The proof will be given in Appendix. 

\begin{lem}\label{lem:3.9}
  Let $e\in V$ be an Ising vector.
  Let $b\in V_2$ and $y \in V_N$ be highest weight vectors for
  $\vir(e)$ such that $2e_{(1)}b=b$ and $e_{(1)}y=hy$, $h\in \{ 0,1/2\}$.
  Then for any $k\in \Z$, there exist $P_{1/2-h,k,j} \in \mathrm{U}(\vir(e)_-)[j]$
  such that
  \[
    \l( 2e_{(i)}-\delta_{i,1}(1-2h)\r)
    \l(b_{(k)}y +\sum_{j>0} P_{1/2-h,k,j}\, b_{(k+j)} y\r)=0
  \]
  for $i\geq 1$.
  Moreover, $P_{1/2-h,k,j}$ are determined only by $h$, $k$, $j$ and $b_{(N+1)}y\in V_0$.
\end{lem}

Suppose $e$ is an Ising vector of a VOA $V$ of OZ-type.
An Ising vector $e$ is said to be \textit{of $\sigma$-type on $V$} if there exists no irreducible
$\la e\ra$-submodule of $V$ isomorphic to $L(\shf,\sfr{1}{16})$.
In this case, we have an isotypical decomposition
\begin{equation}\label{eq:3.3}
  V= V[0]_e\oplus V[\shf]_e
\end{equation}
where $V[h]_e$ is the sum of all irreducible $\la e\ra$-submodules isomorphic to $L(\shf,h)$,
and the Griess algebra of $V$ decomposes into a direct sum of eigenspaces
\begin{equation}\label{eq:3.4}
  V_2=\C e \oplus V_2[0]_e\oplus V_2[\shf]_e
\end{equation}
with $V_2[h]_e=\ker_{V_2}(\o(e)-h)$.
By the fusion rules of $L(\shf,0)$-modules in \eqref{eq:3.2} and
based on the decomposition \eqref{eq:3.3},  we can define an automorphism by
\begin{equation}\label{equation:3.5}
  \sigma_{e}:=(-1)^{2\o(e)} =
  \begin{cases}
    ~1 & \text{on }~ V[0]_e,
    \\
    -1 & \text{on }~ V[\sfr{1}{2}]_e.
  \end{cases}
\end{equation}
The involution $\sigma_e$ is called a \textit{Miyamoto involution} of $\sigma$-type or a \textit{$\sigma$-involution}
(cf.~\cite{Mi}).
By the definition, we have the following conjugation.

\begin{prop}\label{prop:3.2}
  Let $e\in V$ be an Ising vector of $\sigma$-type and $g\in \aut(V)$.
  Then we have $\sigma_{ge}=g\sigma_e g^{-1}$.
\end{prop}

The local structures of subalgebras generated by two Ising vectors of $\sigma$-type are
completely determined in \cite{Mi,Ma}.

\begin{prop}[\cite{Mi,Ma}]\label{prop:3.3}
  Let $V$ be a VOA of OZ-type and let $a$ and $b$ be distinct Ising vectors of
  $\sigma$-type on $V$.
  Then the Griess subalgebra $B$ generated by $a$ and $b$ is one of the following.
  \\
  $(i)$~$(a\mymid b)=0$, $a_{(1)}b=0$ and $B=\C a \oplus \C b$.
  In this case,  $\sigma_a$ commutes with $\sigma_b$ on $V$.
  \\
  $(ii)$~$(a\mymid b)=2^{-5}$, $\sigma_a b=\sigma_b a$, $4a_{(1)}b=a+b-\sigma_a b$ and
  $B=\C a \oplus \C b \oplus \C \sigma_a b$.
  In this case, $\sigma_a\sigma_b$ has order three on $V$.
\end{prop}

\begin{cond}\label{cond:1}
  We consider a VOA $V$ satisfying the following conditions.
  \\
  (1)~$V$ is of OZ-type.
  \\
  (2)~$E_V$ is the set of Ising vectors of $V$ of $\sigma$-type.
  \\
  (3)~$V$ is generated by $E_V$.
\end{cond}

For a VOA $V$ satisfying Condition \ref{cond:1}, we set 
$I_V=\{ \sigma_e \mid e\in E_V\}$ and $G_V=\la I_V\ra$.
It follows from Propositions \ref{prop:3.2} and \ref{prop:3.3} that
$(G_V,I_V)$ is a 3-transposition group.
The 3-transposition groups arising in this manner is classified in \cite{Ma} under the
assumption that $V$ has a compact real form.
They are known to be finite 3-transposition groups of symplectic type (cf.~\cite{CH}).

\begin{rem}\label{rem:3.6}
  If $(G_V,I_V)$ is decomposable then we have a non-trivial partition $I_V=I_1\sqcup I_2$ 
  such that $G_i=\la I_i\ra$ are non-trivial 3-transposition subgroups of $G_V$ and 
  $G_V\cong G_1\times G_2$.
  The partition $I_V=I_1\sqcup I_2$ induces a partition $E_V=E_1\sqcup E_2$ with
  $E_i=\{ e \in E_V \mid \sigma_e\in I_i\}$.
  Let $V^i$ be the subVOA generated by $E_i$, $i=1,2$.
  Then $V\cong V^1\tensor V^2$ and $G_i<\aut(V^i)$ (cf.~\cite{Mi,Ma}).
  Therefore, the study of $V$ satisfying Condition \ref{cond:1} is reduced to the case when
  $G_V$ is indecomposable.
  It is shown in Section 3 of \cite{Ma} that $G_V$ is center-free and the correspondence 
  $E_V\ni e\mapsto \sigma_e\in I_V$ is bijective if $(G_V,I_V)$ is indecomposable and 
  non-trivial.
\end{rem}

%

%

Let $V$ be a VOA satisfying Condition \ref{cond:1} and $(G_V,I_V)$ the associated 
3-transposition group.
We set
\begin{equation}\label{eq:3.6}
  \A= \ds \bigsqcup_{k\geq 0} (E_V\times \Z_{\geq 0})^k
  = \{ ((a_1,n_1),\dots,(a_k,n_k)) \mid k\geq 0,~ a_i\in E_V,~ n_i\geq 0\} .
\end{equation}
We define the weight of $\alpha=((a_1,n_1),\dots,(a_k,n_k))\in \A$ by 
$\wt\, \alpha=n_1+\cds +n_k+k\in \Z_{\geq 0}$ and set
\[
  \A_N= \{ \alpha \in \A \mid \wt\, \alpha =N\} .
\]
Then $\A=\ds \bigsqcup_{N\geq 0} \A_N$.
We also set $\ds \A_{\leq N}=\bigsqcup_{n\leq N} \A_n$.
We define a map $\phi_V : \A \to V$ by
\begin{equation}\label{eq:3.7}
  \phi_V((a_1,n_1),\dots,(a_k,n_k))
  = {a_1}_{(-n_1)}\cds {a_k}_{(-n_k)}\vac .
\end{equation}
It follows that $\phi_V(\alpha)\in V_N$ if $\alpha\in \A$ has weight $N$.

Let $W$ be another VOA satisfying Condition \ref{cond:1} such that the Griess algebras  
of $V$ and $W$ are isomorphic.
Let $\theta: V_2\to W_2$ be an isomorphism of the Griess algebras.
Then $\theta$ defines a bijection between $E_V$ and $E_W$ so that 
we can define the map $\phi_W: \A\to W$ by 
\begin{equation}\label{eq:3.8}
  \phi_W((a_1,n_1),\dots,(a_k,n_k))
  = (\theta a_1)_{(-n_1)}\cds (\theta a_k)_{(-n_k)}\vac .
\end{equation}

\begin{df}\label{df:univ}
Let $V$ be a VOA satisfying Condition \ref{cond:1}.
Let $e\in E_V$, $\alpha\in \A$ and $m>0$.
We say $(e,m,\alpha)$ has a \emph{universal expression} if there exists 
an expression $e_{(m)}\phi_V(\alpha)=\sum_j c_j \phi_V(\beta_j)$ in $\Span\, \phi_V(\A)$
such that for any VOA $W$ satisfying Condition \ref{cond:1} with an isomorphism
$\theta: V_2\to W_2$ of the Griess algebras, we have the same expression
$(\theta e)_{(m)} \phi_W(\alpha)=\sum_j c_j \phi_W(\beta_j)$ in $\Span\, \phi_W(\A)$. 
For convenience, we often say that $e_{(m)}\phi_V(\alpha)$ has a universal expression also.   
\end{df}

\begin{lem}\label{lem:3.10}
  Let $V$ be a VOA satisfying Condition \ref{cond:1}.
  Then $(e,m,\alpha)$ has a universal expression for every 
  $e\in E_V$, $\alpha\in \A$ and $m>0$.
\end{lem}

\pf
We prove the lemma by induction on the weight of $\alpha\in \A$.
Since $\Span\, \phi_V(\A_{\leq 2})=\C \vac \oplus \Span\, E_V$, the lemma holds if $N\leq 2$.
Suppose that the lemma holds for $\alpha\in \A_{\leq N}$ with $N\geq 2$.
Suppose $e\in E_V$, $\alpha=((a_1,n_1),\dots,(a_k,n_k)) \in \A_{N+1}$ and $m>0$.
Let us show that $e_{(m)}\phi_V(\alpha)$ has a universal expression.
Set $b=a_1\in E_V$, $n=n_1$, $\alpha'=((a_2,n_2),\dots,(a_k,n_k))\in \A_{N-n}$ 
and $y=\phi_V(\alpha')\in V_{N-n}$. 
Then $\phi_V(\alpha) = b_{(-n)}y$.
We have
\[
  e_{(m)}\phi_V(\alpha)=e_{(m)}b_{(-n)}y=[e_{(m)},b_{(-n)}]y+b_{(-n)}e_{(m)}y.
\]
By induction, $e_{(m)}y=e_{(m)}\phi_V(\alpha')$ has a universal expression.
On the other hand,
\[
  [e_{(m)},b_{(-n)}]y
  =\l( e_{(0)}b\r)_{(m-n)}y+m\l( e_{(1)}b\r)_{(m-n-1)}y
  +\binom{m}{3}\l( e_{(3)}b\r)_{(m-n-3)}y
\]
and again by induction the latter two terms have universal expressions since
$e_{(1)}b\in \Span\, E_V$ and $e_{(3)}b=(e\mymid b)\vac$.
Therefore, it is enough to show that either $e_{(m)}b_{(-n)}y$ or $\l( e_{(0)}b\r)_{(m-n)}y$ 
can be written as a sum of universal expressions.

\textbf{Case 1: $n>0$.}~
In this case, $\alpha'\in \A_{\leq N-1}$ and $y=\phi_V(\alpha')$ has a weight less than $N$.
Then $\alpha''=((e,0),\alpha')\in \A_{\leq N}$ and both $y=\phi_V(\alpha')$ and 
$e_{(0)}y=\phi_V(\alpha'')$ are in $\Span\,\phi_V(\A_{\leq N})$.
Since the $0$-th operator is a derivation, we have 
$\l( e_{(0)}b\r)_{(m-n)}y=e_{(0)}b_{(m-n)}y-b_{(m-n)}e_{(0)}y$.
If $m-n\leq 0$, the right hand side can be expressed by standard expressions in 
$\Span\,\phi_V(\A)$ as in \eqref{eq:3.7}.
If $m-n>0$, we can apply the induction to $(b,\alpha')$ and $(b,\alpha'')$ to obtain that both  
$e_{(0)}b_{(m-n)}y$ and $b_{(m-n)}e_{(0)}y$ have universal expressions.
Therefore, $\l( e_{(0)}b\r)_{(m-n)}y$ has a universal expression.

\textbf{Case 2: $n=0$ and $m>1$.}~
In this case,  $y\in \phi_V(\A_N)$ and
\[
\begin{array}{ll}
  \l( e_{(0)}b\r)_{(m-n)}y
  &=\l( e_{(0)}b\r)_{(m)}y
  =[e_{(1)},b_{(m-1)}]y-\l( e_{(1)}b\r)_{(m-1)}y
  \medskip\\
  &= e_{(1)}b_{(m-1)}y-b_{(m-1)}e_{(1)}y-\l( e_{(1)}b\r)_{(m-1)}y.
\end{array}
\]
Since $m>1$, we can apply the induction to
$b_{(m-1)}y$, $e_{(1)}y$ and $\l( e_{(1)}b\r)_{(m-1)}y$, and
all of them have universal expressions in $\Span\,\phi_V(\A_{\leq N})$,
and so are $e_{(1)}b_{(m-1)}y$ and $b_{(m-1)}e_{(1)}y$ again by induction.

\textbf{Case 3: $n=0$ and $m=1$.}~
By linearity, it is enough to show that either $e_{(m)}b_{(-n)}y$ or $\l( e_{(0)}b\r)_{(m-n)}y$ 
can be written as a sum of universal expressions when 
$b\in \Span\,E_V$ and $y\in \Span\, \phi_V(\A_N)$.
Note that the choice of $b$ and $y$ does not matter in showing that the expressions 
to be obtained below are universal.
We can also assume that $e_{(1)}b=t b$ with $t=0$, 1/2 or 2.
If $t=2$ then $b\in \C e$ and the claim follows since $(e_{(0)}e)_{(1)}=-e_{(0)}$.
If $t=0$ then $b$ is a highest weight vector of highest weight 0 so that
$e_{(0)}b=0$ and the claim follows.
Therefore we assume that $b$ is a highest weight vector for $\vir(e)$
with highest weight 1/2.
In this case, we consider $e_{(1)}b_{(0)}y$ instead of $(e_{(0)}b)_{(1)}y$.
Since $\la e\ra\cong L(\shf,0)$ is rational, $e_{(1)}$ acts on $\phi_V(\A_{\leq N})$
semisimply and $y$ is a sum of descendants $e_{(-s_1)}\cds e_{(-s_k)}v_\lambda$,
$s_1\geq \cds \geq s_k \geq 0$, $k\geq 0$, of highest weight vectors
$v_\lambda\in \Span\,\phi_V(\A_{\leq N})$ for $\vir(e)$.
By the inductive assumption, this expression is universal, and again by linearity
we may assume that $y=e_{(-s_1)}\cds e_{(-s_k)}v_\lambda$ with
$s_1\geq \cds \geq s_k\geq 0$ and $k\geq 0$.
Moreover, by Lemma \ref{lem:3.8}, we may assume that $s_{k-1}\geq 1$ when $s_k=0$.
In other words, $s_1=0$ implies $k=1$.
Since $e$ is of $\sigma$-type, the highest weight of $v_\lambda$ is either 0 or 1/2.
We further divide Case 3 into three subcases.

\textbf{Case 3-1: $k=0$.}~
Suppose $y=v_\lambda\in \Span\, \phi_V(\A_N)\subset V_N$ is a highest weight vector for $\vir(e)$.
By the inductive assumption,  $b_{(N+1)}y\in V_0$ has a universal expression and is
uniquely determined by the Griess algebra of $V$.
Then it follows from Lemma \ref{lem:3.9} and the inductive assumption that
$e_{(1)}b_{(0)}y$ has a universal expression.

\textbf{Case 3-2: $k>0$ and $s_1>0$.}~
Suppose $k>0$. 
Then $y=e_{(-s)}w$ with $s=s_1\geq 0$ and 
$w=e_{(-s_2)}\cds e_{(-s_k)}v_\lambda\in \Span\,\phi_V(\A_{N-s-1})$.
Suppose $s>0$.
We have
\[
  e_{(1)}b_{(0)}y
  = e_{(1)}b_{(0)}e_{(-s)}w
  = e_{(1)}[b_{(0)},e_{(-s)}]w +e_{(1)}e_{(-s)}b_{(0)}w.
\]
We also have 
\[
  e_{(1)}e_{(-s)}b_{(0)}w=e_{(-s)}e_{(1)}b_{(0)}w+[e_{(1)},e_{(-s)}]b_{(0)}w
  =e_{(-s)}e_{(1)}b_{(0)}w+(s+1)e_{(-s)}b_{(0)}w .
\]
Since $b_{(0)}w\in \Span\,\phi_V(\A_{N-s})$, $e_{(1)}b_{(0)}w$ has a universal
expression by induction.
Therefore $e_{(1)}e_{(-s)}b_{(0)}w$ has a universal expression.
On the other hand, we have
\[
\begin{array}{l}
  e_{(1)}[b_{(0)},e_{(-s)}]w
  =-e_{(1)}[e_{(-s)},b_{(0)}]w
  \medskip\\
  \ds =-e_{(1)}\l(\l( e_{(0)}b\r)_{(-s)}-s\l( e_{(1)}b\r)_{(-s-1)}
   +\binom{-s}{3}\l( e_{(3)}b\r)_{(-s-3)}\r) w
  \medskip\\
  \ds = -e_{(1)}e_{(0)}b_{(-s)}w+e_{(1)}b_{(-s)}e_{(0)}w+se_{(1)}\l( e_{(1)}b\r)_{(-s-1)}w
  -\binom{-s}{3}e_{(1)}(e\mymid b)\vac_{(-s-3)}w
  \medskip\\
  =-e_{(0)}b_{(-s)}w-e_{(0)}e_{(1)}b_{(-s)}w
  +e_{(1)}b_{(-s)}e_{(0)}w+se_{(1)}\l( e_{(1)}b\r)_{(-s-1)}w.
\end{array}
\]
By the inductive assumption, $e_{(1)}b_{(-s)}w$ has a universal expression,
and by Case 1, the terms $e_{(1)}b_{(-s)}e_{(0)}w$ and $e_{(1)}\l( e_{(1)}b\r)_{(-s-1)}w$ also have
universal expressions.
Therefore, all the terms in the above have universal expressions and the claim follows.

\textbf{Case 3-3: $k=1$ and $s_1=0$.}~
Now it remains to show that $e_{(1)}b_{(0)}y$ has a universal expression when
$y=e_{(0)}v_\lambda$.
Note that $v_\lambda\in \Span\,\phi_V(\A_{N-1})$.
We have
\[
\begin{array}{ll}
  e_{(1)}b_{(0)}e_{(0)}v_\lambda
  &= e_{(1)}e_{(0)}b_{(0)}v_\lambda + e_{(1)}[b_{(0)},e_{(0)}]v_\lambda
  \medskip\\
  &= e_{(0)}b_{(0)}v_\lambda +e_{(0)}e_{(1)}b_{(0)}v_\lambda -e_{(1)}[e_{(0)},b_{(0)}]v_\lambda
  \medskip\\
  &= e_{(0)}b_{(0)}v_\lambda +e_{(0)}e_{(1)}b_{(0)}v_\lambda
  -e_{(1)}\l(e_{(0)}b\r)_{(0)}v_\lambda .
\end{array}
\]
Since $e_{(0)}b_{(0)}v_\lambda \in \Span\,\phi_V(\A_{N+1})$
and $e_{(1)}b_{(0)}v_\lambda$ has a universal expression by the inductive assumption,
it suffices to consider the term $e_{(1)}\l(e_{(0)}b\r)_{(0)}v_\lambda$.
We have
\[
\begin{array}{l}
  e_{(1)}\l(e_{(0)}b\r)_{(0)}v_\lambda
  = \l(e_{(0)}b\r)_{(0)}e_{(1)}v_\lambda+\l[ e_{(1)},\l(e_{(0)}b\r)_{(0)}\r] v_\lambda
  \medskip\\
  = \l( e_{(0)}b_{(0)}-b_{(0)}e_{(0)}\r) e_{(1)}v_\lambda
  +\l({e_{(0)}}^2\, b\r)_{(1)} v_\lambda +\l(e_{(1)}e_{(0)}b\r)_{(0)} v_\lambda
  \medskip\\
  = e_{(0)}b_{(0)}e_{(1)}v_\lambda -b_{(0)}e_{(0)}e_{(1)}v_\lambda
  +\l({e_{(0)}}^2b\, \r)_{(1)} v_\lambda +\dfr{3}{2}\l(e_{(0)}b\r)_{(0)} v_\lambda
  \medskip\\
  = e_{(0)}b_{(0)}e_{(1)}v_\lambda -b_{(0)}e_{(0)}e_{(1)}v_\lambda
  +\l({e_{(0)}}^2b\, \r)_{(1)} v_\lambda +\dfr{3}{2}\l( e_{(0)}b_{(0)}
  -b_{(0)}e_{(0)}\r) v_\lambda .
\end{array}
\]
By induction, all the terms above except $\l({e_{(0)}}^2\,b\r)_{(1)} v_\lambda$
have universal expressions.
Since $3{e_{(0)}}^2b=4e_{(-1)}b$ by Lemma \ref{lem:3.8},
it is enough to rewrite the term $\l(e_{(-1)}b\r)_{(1)} v_\lambda$.
By the iterate formula, we have
\[
\begin{array}{l}
  \l(e_{(-1)}b\r)_{(1)} v_\lambda
  = \dsum_{i\geq 0}\l( e_{(-1-i)}b_{(1+i)}+b_{(-i)}e_{(i)}\r) v_\lambda
\end{array}
\]
and by induction all the terms in the right hand side have universal expressions.
This completes the proof of Lemma \ref{lem:3.10}.
\qed

The following proposition is proved in Lemma 3.5 of \cite{JL} in the case 
when $G_V$ is a symmetric group.

\begin{prop}[\cite{JL}]\label{prop:3.5}
  Let $V$ be a VOA satisfying Condition \ref{cond:1}.
  Then $V$ is linearly spanned by
  \[
    \{ {a_1}_{(-n_1)}\cds {a_k}_{(-n_k)}\vac
    \mid k\geq 0,~ a_i \in E_V,~  n_i\geq 0\} ,
  \]
  i.e., $V=\Span\, \phi_V(\A)$.
\end{prop}

\pf
It follows from Lemma \ref{lem:3.10} and the Borcherds identity (cf.~\cite{MaN}) 
that the subspace $\Span\,\phi_V(\A)$ forms a subalgebra of $V$.
Since $E_V \subset \phi_V(\A_2)$ and $V$ is generated by $E_V$, we obtain $V=\Span\,\phi_V(\A)$.
\qed

\begin{prop}\label{prop:3.10}
  If $V$ is a simple VOA satisfying Condition \ref{cond:1}, then its Griess algebra 
  is isomorphic to the non-degenerate quotient of the Matsuo algebra 
  $B_{\sshf,\sshf}(G_V)$ associated with $G_V=\la I_V\ra$, $I_V=\{ \sigma_e \mid e\in E_V\}$.
\end{prop}

\pf
It follows from Propositions \ref{prop:3.3} and \ref{prop:3.5} that the mapping 
$\pi_V(x^i)=e$, where $i=\sigma_e$ in $I_V$, defines an epimorphism
\[
  \pi_V : B_{\sshf,\sshf}(G_V)\to V_2=\Span\,E_V
\]
of commutative algebras which preserves the bilinear forms.
Since $V$ is simple, the bilinear form on $V_2$ is non-degenerate so that the kernel of
$\pi_V$ is exactly the radical of the bilinear form on $B_{\sshf,\sshf}(G_V)$ (cf.~\cite{Ma}).
Therefore, the Griess algebra of $V$ is isomorphic to the non-degenerate quotient of 
$B_{\sshf,\sshf}(G_V)$.
\qed

Let $V^1$ and $V^2$ be simple VOAs satisfying Condition \ref{cond:1}.
Assume that the Griess algebras of $V^1$ and $V^2$ are isomorphic.
We shall prove that $V^1$ and $V^2$ are isomorphic.
Let $E_i$ be the set of Ising vectors of $V^i$ of $\sigma$-type.
Since the Griess algebras of $V^1$ and $V^2$ are isomorphic,
there exists a bijection $\theta:E_1\to E_2$.
We identify $E_1$ with $E_2$ via $\theta$ and consider the index set $\A$ associated with 
$E_1=E_2$ in \eqref{eq:3.6}.

\begin{prop}\label{prop:3.7}
  If $V^1$ and $V^2$ are VOAs satisfying Condition \ref{cond:1} with the same Griess algebra, 
  then $(\phi_{V^1}(\alpha)\mymid \phi_{V^1}(\beta))=(\phi_{V^2}(\alpha)\mymid \phi_{V^2}(\beta))$ 
  for all $\alpha$, $\beta\in \A$.
\end{prop}

\pf
By the invariance
$(e_{(-n)}\phi_{V^i}(\alpha) \mymid \phi_{V^i}(\beta)) = 
(\phi_{V^i}(\alpha) \mymid e_{(n+2)}\phi_{V^i}(\beta))$ for 
$e\in E_{V^i}$ and $\alpha$, $\beta\in \A$, 
the proposition follows from Lemma \ref{lem:3.10} and the definition of 
the universal expression.
\qed

Since $V^i$ is of OZ-type, $V^i$ has a unique invariant bilinear form such that $(\vac|\vac)=1$, 
which is non-degenerate as $V^i$ is simple.
We define a linear map $f:V^1\to V^2$ by
\[
  f(\phi_{V^1}(\alpha))=\phi_{V^2}(\alpha)~~~\mbox{for $\alpha\in \A$}.
\]

\begin{lem}\label{lem:3.13}
  $f$ is a well-defined linear isomorphism.
\end{lem}

\pf
We extend the maps $\phi_{V^1}$, $\phi_{V^2}$ to linear maps from the span of $\A$.
Suppose $\alpha\in \ker\, \phi_{V^1}$ and let $\beta\in \A$.
Then $(\phi_{V^1}(\alpha)\mymid \phi_{V^1}(\beta))=0$.
By Proposition \ref{prop:3.7}, we have $(\phi_{V^2}(\alpha)\mymid \phi_{V^2}(\beta))=0$.
As it holds for all $\beta\in \A$, and $\phi_{V^2}(\beta)$ runs over the spanning set of $V^2$ 
by Proposition \ref{prop:3.5}, we have $\phi_{V^2}(\alpha)=0$ since the form on $V^2$ 
is non-degenerate.
Therefore, $f$ induces a linear isomorphism between 
$V^1\cong \Span_{V^1}\, \A/\ker\, \phi_{V^1}$ 
and $V^2\cong \Span_{V^2}\, \A/\ker\, \phi_{V^2}$.
\qed

Note that by restriction $f$ defines an isomorphism between Griess algebras of $V^1$ and $V^2$.
It follows from Proposition 5.7.9 of \cite{LL} that $f$ is a homomorphism of VOAs if 
and only if $f(e_{(m)}v)=f(e)_{(m)}f(v)$ for $e\in E_{V^1}$, $v\in V^1$ and $m\in \Z$.
This is clear if $m\leq 0$, and also true for $m>0$ by Lemma \ref{lem:3.10} and 
the definition of the universal expression.
Thus, we have proved the following theorem.

\begin{thm}\label{thm:3.6}
  If $V^1$ and $V^2$ are simple VOAs of OZ-type satisfying Condition \ref{cond:1} with 
  isomorphic Griess algebras, then the isomorphism extends to an isomorphism of VOAs. 
\end{thm}

As a by-product of the theorem above, we give a slight simplification of Matsuo's 
classification of certain 3-transposition groups realizable by VOAs.

\begin{prop}\label{prop:3.12}
  Let $V$ be a VOA satisfying Condition \ref{cond:1} and
  let $a$ and $b$ be Ising vectors such that $(a\mymid b)=2^{-5}$.
  Then the subalgebra $\la a,b\ra$ is simple and isomorphic to
  \[
    L(\shf,0)\tensor L(\sfr{7}{10},0)
    \oplus L(\shf,\shf)\tensor L(\sfr{7}{10},\sfr{3}{2}).
  \]
\end{prop}

\pf
By Proposition \ref{prop:3.3} the Griess algebra of $\la a,b\ra$ is
isomorphic to the Matsuo algebra associated with the symmetric group
$\mathfrak{S}_3$ of degree three.
Since $\la a,b\ra$ is of OZ-type, and the radical $J$ of the invariant bilinear
form on $\la a,b\ra$ is the unique maximal ideal, it suffices to show that $J=0$.
The subalgebra $\la a,b\ra$ has a conformal vector
\begin{equation}\label{eq:3.8+}
  \eta =\dfr{4}{5}(a+b+\sigma_a b)
\end{equation}
given by \eqref{eq:2.3}.
Since $V$ is of OZ-type,  it follows from Theorem 5.1 of \cite{FZ} that 
the conformal vector of $V$ is an orthogonal sum of $\eta$ and its complement.
Therefore, we have a grading $J=\oplus_{n\geq 0} J_n$ with $J_n=J\cap V_n$.
Let $h$ be the top weight of $J$, i.e., $J_h\ne 0$ and $J_n=0$ for $n<h$.
Since $a$ and $b$ are of $\sigma$-type, the zero-modes  $\o(a)$,
$\o(b)$ and $\o(\sigma_a b)$ act on $J_h$ semisimply with eigenvalues in
$\{ 0,1/2\}$.
Therefore, we have
\begin{equation}\label{eq:3.9}
  \tr_{J_h}\, \o(\eta)= \dfr{4}{5} \tr_{J_h}\,\o(a+b+\sigma_ab)
  \leq \dfr{4}{5}\cd 3\cd \dfr{1}{2}\dim J_h =\dfr{6}{5}\dim J_h.
\end{equation}
On the other hand, since $\eta$ is the conformal vector of $\la a,b\ra$,
the zero mode $\o(\eta)$ acts by the top weight $h$ on $J_h$.
Therefore $\tr_{J_h}\,\o(\eta)=h\dim J_h$ and from \eqref{eq:3.9}
we obtain $h\leq 6/5$, showing $h=1$.
This contradicts the assumption that $V$ is of OZ-type.
Therefore $J=0$ and $\la a,b\ra$ is simple.
Then by Theorem \ref{thm:3.6} its VOA structure is uniquely determined as
in the assertion (cf.~\cite{LY}).
\qed

\begin{prop}[\cite{Ma}]\label{prop:3.13}
  Let $V$ be a VOA satisfying Condition \ref{cond:1}.
  Let $E_V$ be the set of Ising vectors of $V$ of $\sigma$-type and
  set $G_V=\la \sigma_e \mymid e\in E_V\ra$.
  Then $G_V$ is a 3-transposition group of symplectic type.
\end{prop}

\pf
It is shown in Proposition 1 of \cite{Ma}\footnote{Proposition 3.3.8 of
\href{http://arxiv.org/abs/math/0311400}{\tt arXiv:math/0311400}.}
that $G_V$ is of symplectic type provided that $V$ has a compact real form
containing $E_V$.
The key idea in the proof of \cite{Ma} is to find a non-zero highest weight
vector for $L(\sfr{7}{10},0)$ with highest weight $7/10$ when $G$ is not
of symplectic type.
In (loc.~cit.), the compact real form is used only to show the existence
of the subalgebra isomorphic to $L(\sfr{7}{10},0)$, and the existence of
the highest weight vector depends only on the structure of the Griess algebra
and this part is independent of the compact real form.
Now by Proposition \ref{prop:3.12}, we can obtain $L(\sfr{7}{10},0)$ without
the assumption on the compact real form, and we obtain the same contradiction
as in \cite{Ma} if $G_V$ is not of symplectic type.
\qed

\begin{cor}\label{cor:3.14}
  Let $V$ be a simple VOA satisfying Condition \ref{cond:1} and
  let $V_\R$ be the real VOA generated by the set $E_V$ of Ising vectors
  of $V$ of $\sigma$-type.
  If the non-degenerate quotient of the real Matsuo algebra $B_{\sshf,\sshf}(G_V)_\R$
  associated with $G_V=\la \sigma_e \mymid e\in E_V\ra$ is positive definite,
  then $V_\R$ is a compact real form of $V$.
  In this case a non-trivial indecomposable component of the 3-transposition
  group $G_V$ is isomorphic to one of the groups listed in Theorem 1 of \cite{Ma}.
\end{cor}

\pf
By the assumption, the invariant bilinear form of the Griess algebra of $V_\R$ is
positive definite so that $E_V$ is a finite set since $(a\mymid b)=2^{-5}$ or $0$
for distinct $a$, $b\in E_V$.
It follows from Proposition \ref{prop:3.13} and $V=\la E_V\ra$ that
$G_V$ is a finite 3-transposition group of symplectic type.
By the positivity of the Griess algebra, it follows from \cite{Ma} that
a non-trivial indecomposable component of $G_V$ is one of the groups listed
in Theorem 1 of \cite{Ma}.
Since all the examples of VOAs in (loc.~cit.) have compact real forms generated
by Ising vectors, it follows from the uniqueness of the VOA structure of $V$
shown in Theorem \ref{thm:3.6} that $V_\R$ is a compact real form of $V$.
\qed

\begin{rem}
  In Theorem 1 of \cite{Ma}, there exist non-isomorphic 3-transposition groups
  but realized by the same VOA.
  For example, $\mathrm{O}_{10}^+(2)$ and $2^8{:}\mathrm{O}_8^+(2)$ are realized
  by $V_{\sqrt{2}E_8}^+$.
  This is because $2^8{:}\mathrm{O}_8^+(2)<\mathrm{O}_{10}^+(2)$ and
  $B_{\sshf,\sshf}(2^8{:}\mathrm{O}_8^+(2))$ is a subalgebra of
  $B_{\sshf,\sshf}(\mathrm{O}_{10}^+(2))$ but they have isomorphic non-degenerate quotients.
  In our argument, we always take $E_V$ to be the set of \emph{all} Ising vectors
  of $V$ of $\sigma$-type, and for $V=V_{\sqrt{2}E_8}^+$ we obtain the maximal one
  $G_V=\mathrm{O}_{10}^+(2)$.
\end{rem}

We propose a conjecture on positivity of a simple VOA satisfying Condition \ref{cond:1}.

\begin{conj}\label{conj:3.16}
  Let $V$ be a simple VOA satisfying Condition \ref{cond:1}.
  Then the bilinear form on the $\R$-span of $E_V$ is positive definite,
  i.e., the non-degenerate quotient of the real Matsuo algebra
  $B_{\sshf,\sshf}(G_V)_\R$ associated with $G_V$ is positive definite.
\end{conj}

If this conjecture is true, it follows from Corollary \ref{cor:3.14} that
the classification of VOA-realizable 3-transposition groups together with VOAs
generated by Ising vectors of $\sigma$-type in \cite{Ma} is complete without
the assumption on the compact real form of a VOA.

\begin{rem}
  It follows from \cite{R} that for any 3-transposition group $G$ there exists a simple 
  VOA of OZ-type whose weight two subspace contains a quotient of the Matsuo 
  algebra $B_{\sshf,\sshf}(G)$ as a subalgebra of the Griess algebra.
  It follows from Theorem \ref{thm:3.6} that if $V$ has a compact real form and
  $G$ is not in the list of Matsuo's classification \cite{Ma} then such a VOA does not satisfy 
  Condition \ref{cond:1}.
\end{rem}

\section{Simplicity of type $A_n$}

Let $\Phi(A_n)$ be the root system of type $A_n$.
We fix a system of simple roots $\alpha_1,\dots,\alpha_n$ of $\Phi(A_n)$ such that
$(\alpha_i\mymid \alpha_j)=-1$ if and only if $\abs{i-j}=1$, and denote by $\Phi(A_n)^+$
the set of positive roots.
We also fix root subsystems $\Phi(A_i)$ for $1\leq i\leq n$ to be the root systems generated
by $\{ \alpha_j \mymid 1\leq j\leq i\}$ and set $\Phi(A_i)^+=\Phi(A_n)^+\cap \Phi(A_i)$.
Let $r_\alpha$ be the reflection associated with $\alpha \in \Phi(A_n)$.
The Weyl group $W(A_i)$ of $\Phi(A_i)$ is a 3-transposition group isomorphic to
the symmetric group $\sym_{i+1}$ of degree $i+1$ with the set of transpositions
$\{ r_\alpha\mid \alpha\in \Phi(A_i)^+\}$.

Let $M_{A_n}$ be the commutant subalgebra of the diagonal subalgebra
$L_{\widehat{\mathfrak{sl}}_2}(n+1,0)$ in the tensor product
$L_{\widehat{\mathfrak{sl}}_2}(1,0)^{\tensor n+1}$ (cf.~\cite{LS,JL}).
Since $L_{\widehat{\mathfrak{sl}}_2}(1,0)$ has a compact real form, so does
$L_{\widehat{\mathfrak{sl}}_2}(1,0)^{\tensor n+1}$ and hence $M_{A_n}$ is a simple VOA with a compact 
real form.
It is known (cf.~\cite{DLMN,Ma,LSY,JL}) that if $n>1$ then $M_{A_n}$ is a rational VOA satisfying
Condition \ref{cond:1} with $G_{M_{A_n}}=W(A_n)=\sym_{n+1}$ and the Griess algebra of $M_{A_n}$
is isomorphic to the Matsuo algebra $B_{\sshf,\sshf}(\sym_{n+1})$.
In particular, $B_{\sshf,\sshf}(\sym_{n+1})$ is non-degenerate.
The set of Ising vectors of $M_{A_n}$ is in one-to-one correspondence with  the set of
reflections $\{ r_\alpha \mymid \alpha \in\Phi(A_n)^+\}$ so that we can index the set
of Ising vectors of $M_{A_n}$ as $E_{M_{A_n}}=\{ x^\alpha \mymid \alpha\in \Phi(A_n)^+\}$.
Note that $\sigma_{x^\alpha}=r_\alpha$ in $\aut(M_{A_n})=W(A_n)$ if $n>1$.

Let $V$ be a VOA satisfying Condition \ref{cond:1} such that the
associated 3-transposition group $G_V=\la \sigma_e\mymid e\in E_V\ra$ is
isomorphic to $W(A_n)\cong \sym_{n+1}$.
In this case the Griess algebra of $V$ is isomorphic to the non-degenerate Matsuo algebra
$B_{\sshf,\sshf}(\sym_{n+1})$ and the set of Ising vectors of $V$ is given by
$E_V=\{ x^\alpha \mymid \alpha \in \Phi(A_n)^+\}$.
By Theorem \ref{thm:3.6}, $V$ has the unique simple quotient isomorphic to $M_{A_n}$.
The purpose of this section is to prove the following theorem.

\begin{thm}\label{thm:4.6}
  Let $V$ be a VOA satisfying Condition \ref{cond:1} such that
  $G_V\cong \mathfrak{S}_{n+1}$.
  Then $V$ is simple and isomorphic to $M_{A_n}$.
\end{thm}


Set $V^{[i]}:=\la x^\alpha \mymid \alpha \in \Phi(A_i)^+\ra$ for $1\leq i\leq n$.
Then we have a tower of subVOAs
\begin{equation}\label{eq:4.1}
  V^{[1]}\subset V^{[2]}\subset \cds \subset V^{[n]}=V
\end{equation}
such that $V^{[i]}$ satisfies Condition \ref{cond:1} with $G_{V^{[i]}}= W(A_i)=\sym_{i+1}$ for $i\geq 2$.
We know $V^{[1]}\cong L(\shf,0)$ is simple and it follows from Proposition \ref{prop:3.12}
that $V^{[2]}$ is also simple.
Let $\w^i$ be the conformal vector of $V^{[i]}$ given by \eqref{eq:2.3}.
We also set $\w^0=0$ for convention.
Then $\eta^i:= \w^i-\w^{i-1}$ for $1\leq i\leq n$ are mutually orthogonal Virasoro vectors
with central charges $c_i$ in \eqref{eq:3.1} (cf.~\cite{DLMN,Ma,FZ}).

Next we consider the Zhu algebra $A(V)=V/O(V)$ of $V$ (cf.~\cite{Z}). 
We study a natural homomorphism from $A(V)$ to $A(M_{A_n})$ and determine possible 
irreducible subquotients of the adjoint module of $V$.
In particular, we show that $V$ has a composition series.
Our argument depends heavily on the explicit structure of the Zhu algebra $A(M_{A_n})$ 
determined in \cite{JL}.
We denote by $[a]$ the residue class $a+O(V)$ of $a\in V$ in $A(V)$.
The following proposition is shown in Theorem 3.6 of \cite{JL} in the case when $G_V$ is 
isomorphic to a symmetric group but based on  Proposition \ref{prop:3.5} the same proof 
works fine for any VOA satisfying Condition \ref{cond:1}.

\begin{lem}[\cite{JL}]\label{lem:4.1}
  If $V$ is a VOA satisfying Condition \ref{cond:1} then the Zhu algebra of $V$ 
  is generated by $E_V+O(V)$.
\end{lem}

We say a $V$-module $N$ is \emph{of $\sigma$-type}\footnote{Such a module is called a  
type I module in \cite{JL}.} if for each $x^\alpha\in E_V$, there is no 
$\la x^\alpha\ra$-submodule of $N$ isomorphic to $L(\shf,\sfr{1}{16})$.
Let $\mathfrak{a}$ be the two-sided ideal of $A(V)$ generated by
$[x^\alpha][2x^\alpha-\vac]$ with $\alpha \in \Phi(A_n)^+$.
If $N=\oplus_{k\geq 0}N(k)$ is an $\N$-gradable $V$-module of $\sigma$-type,
then the top level $N(0)$ is an $A(V)$-module on which $[x^\alpha]$ with $\alpha \in \Phi(A_n)^+$
acts semisimply with possible eigenvalues in $\{ 0,1/2\}$ so that the ideal $\mathfrak{a}$
acts trivially on it and $N(0)$ is a module over the quotient $A(V)/\mathfrak{a}$.
Since $M_{A_n}$ is the simple quotient of $V$ by Theorem \ref{thm:3.6},
$A(M_{A_n})$ is a quotient algebra of $A(V)$ and we can define $A(M_{A_n})/\mathfrak{a}$ as well.

Let $\{ r_\alpha \mymid \alpha\in \Phi(A_n)^+\}$ be the set of transpositions of
$W(A_n)=\sym_{n+1}$ and let $T_{n+1}$ be the quotient algebra of $\C[\sym_{n+1}]$ by the ideal
$\mathfrak{b}$ generated by
\[
  \{ (r_\alpha+r_\beta+r_\alpha r_\beta r_\alpha)(r_\alpha+r_\beta+r_\alpha r_\beta r_\alpha-3) \mid
  \alpha, \beta\in \Phi(A_n)^+,~ r_\alpha \sim r_\beta\} .
\]
It is shown in \cite{JL} that $A(M_{A_n})/\mathfrak{a}$ is isomorphic to $T_{n+1}$.

\begin{lem}[\cite{JL}]\label{lem:4.2}
  The natural map $A(V)\to A(M_{A_n})$ induces an isomorphism 
  between $A(V)/\mathfrak{a}$ and $A(M_{A_n})/\mathfrak{a}$.
\end{lem}

\pf
By the defining relation of the ideal $\mathfrak{a}$, we have
$[\vac-4x^\alpha]^2=[\vac]$ in $A(V)/\mathfrak{a}$.
Let $x^\alpha$, $x^\beta\in E_V$ be Ising vectors such that $r_\alpha \sim r_\beta$.
Then $\la x^\alpha,x^\beta\ra$ is isomorphic to $M_{A_2}$ by Proposition \ref{prop:3.12}.
We denote the Ising vector corresponding to
$r_\alpha \circ r_\beta=r_\alpha r_\beta r_\alpha\in W(A_n)$ by $x^{\alpha\circ \beta}$.
It is known (cf.~\cite{LY}) that $M_{A_2}$ has two non-isomorphic irreducible modules
of $\sigma$-type with top weights 0 and 3/5.
Since the conformal vector of $\la x^\alpha,x^\beta\ra$ is given by
$(4/5)(x^\alpha+x^\beta+x^{\alpha\circ \beta})$ (cf.~Eq.~\eqref{eq:2.3}), it follows that
\[
  [x^\alpha+x^\beta+x^{\alpha\circ \beta}][4(x^\alpha+x^\beta+x^{\alpha\circ \beta})-3\cd \vac]=0
\]
in $A(V)/\mathfrak{a}$.
Let $\mathfrak{c}$ be the ideal generated by $[x^\gamma][2x^\gamma-\vac]$ with 
$\gamma\in \{\alpha, \beta, \alpha\circ \beta\}$ in $A(\la x^\alpha,x^\beta\ra)$.
It is shown in \cite{JL} that
\[
  [\vac -4x^\alpha][\vac -4x^\beta][\vac-4x^\alpha]-[\vac-4x^{\alpha\circ\beta}]=0
\]
in 
$A(\la x^\alpha,x^\beta\ra)/\mathfrak{c}=A(M_{A_2})/\mathfrak{c}$ and hence 
the relation above also holds in $A(V)/\mathfrak{a}$.
It is also shown in Eq.~(3.13) of \cite{JL} that if $r_\alpha r_\beta=r_\beta r_\alpha$
then $[x^\alpha][x^\beta]=[x^\beta][x^\alpha]$, implying the relation
$[\vac -4x^\alpha][\vac -4x^\beta]=[\vac-4x^\beta][\vac-4x^{\alpha}]$.
Therefore, the map $f : r_\alpha \mapsto [\vac -4x^\alpha]$ defines an algebra homomorphism 
from $T_{n+1}=\C[\sym_{n+1}]/\mathfrak{b}$ to $A(V)/\mathfrak{a}$.
Since $A(V)$ is generated by $[x^\alpha]$ with $\alpha\in \Phi(A_n)^+$ by Lemma \ref{lem:4.1},
$f$ is indeed surjective.
Thus we have obtained the following surjections.
\[
  T_{n+1}\xrightarrow{\,~f~} A(V)/\mathfrak{a}
  \xrightarrow{~~~} A(M_{A_n})/\mathfrak{a}.
\]
Since $T_{n+1}$ and $A(M_{A_n})/\mathfrak{a}$ are isomorphic, $f$ is also an isomorphism.
\qed

\begin{cor}\label{cor:4.3}
  An irreducible $V$-module of $\sigma$-type is an irreducible module over the simple
  quotient $M_{A_n}$.
  In particular, there exist finitely many irreducible $V$-modules of $\sigma$-type.
\end{cor}

Since the adjoint $V$-module $V$ is of $\sigma$-type, we have:

\begin{lem}\label{lem:4.4}
  The adjoint $V$-module $V$ has a composition series.
\end{lem}

\begin{lem}\label{lem:4.5}
  The Virasoro vectors $\eta^i$, $1\leq i\leq n$, are simple.
  That is, $\la \eta^i\ra \cong L(c_i,0)$.
\end{lem}

\pf
We prove the lemma by induction on $n$.
It follows from Proposition \ref{prop:3.12} that $\eta^1$ and $\eta^2$ are simple.
Suppose $n>2$ and $V^{[i]}\cong M_{A_i}$ for $1\leq i<n$.
Then $\eta^i$ are simple for $1\leq i<n$.
It suffices to show that $\eta^n$ is simple.
Since $\w^{n-1}=\eta^1+\cds +\eta^{n-1}$ is the conformal vector of $V^{[n-1]}$,
$V^{[n-1]}\tensor \la \eta^n\ra$ is a full subVOA of $V^{[n]}=V$ (cf.~\cite{FZ}).
By induction, $V^{[n-1]}$ is simple and isomorphic to $M_{A_{n-1}}$.
Suppose $\la \eta^n\ra$ is not simple.
Then it follows from the structure of the Verma modules over the Virasoro algebra
(cf.~\cite{IK,W}) that $\la \eta^n\ra$ has a singular vector $u$ of weight $(n+1)(n+2)$.
It follows from Lemma \ref{lem:4.4} that there exists a composition series
$V=J_0 \supset J_1\supset \cds \supset J_k=0$ such that $J_i/J_{i+1}$ is an
irreducible $V$-module of $\sigma$-type.
It follows from Corollary \ref{cor:4.3} that $J_i/J_{i+1}$ are irreducible $M_{A_n}$-modules.
Since $M_{A_n}$ is the unique simple quotient of $V$ by Theorem \ref{thm:3.6},
the top quotient $V/J_1$ is isomorphic to $M_{A_n}$.
It is known (cf.~\cite{LS}) that $\eta^i$ is a simple Virasoro vector in the simple
quotient $M_{A_n}$ so that $u$ is zero in the quotient $V/J_1$, i.e., $u\in J_1$.
So there exists a factor $J_i/J_{i+1}$ such that $u\in J_i$ but $u\not \in J_{i+1}$.
Then $u+J_{i+1}$ is a highest weight vector for the simple Virasoro vector $\eta^n+J_1$ of
$V/J_1\cong M_{A_n}$ with highest weight $(n+1)(n+2)$, which contradicts the classification
of irreducible $L(c_n,0)$-modules (cf.~\eqref{eq:3.1}).
Thus $\eta^n$ is simple in $V$.
\qed
\medskip

Now we prove Theorem \ref{thm:4.6} by induction on $n>2$. 
We assume that $V^{[i]}\cong M_{A_i}$ for $1\leq i<n$.
By the definition, $M_{A_i}\tensor L_{\widehat{\mathfrak{sl}}_2}(i+1,0)$ is a full subVOA of
$L_{\widehat{\mathfrak{sl}}_2}(1,0)^{\tensor i+1}$.
It is shown in \cite{JL} that all the irreducible $M_{A_i}$-modules of $\sigma$-type are
given by the decomposition
\begin{equation}\label{eq:4.2}
  L_{\widehat{\mathfrak{sl}}_2}(1,0)^{\tensor i+1}
  = \bigoplus_{0\leq 2j\leq i+1} M_{A_i}(2j)\tensor L_{\widehat{\mathfrak{sl}}_2}(i+1,2j),
\end{equation}
where $M_{A_i}(2j)$, $0\leq 2j\leq i+1$, are inequivalent irreducible $M_{A_i}$-modules
of $\sigma$-type.
The top weight of $M_{A_i}(2j)$ is positive if $j>0$ and $M_{A_i}(0)$ is isomorphic to
the adjoint module $M_{A_i}$.
By Lemma \ref{lem:4.5}, $V=V^{[n]}$ has a full subVOA $V^{[n-1]}\tensor \la \eta^n \ra$ isomorphic
to $M_{A_{n-1}}\tensor L(c_n,0)$ which is rational by \cite{JL,W}.
It is shown in \cite{GKO,LS} that $M_{A_n}(2j)$ has the following decomposition as
an $M_{A_{n-1}}\tensor L(c_n,0)$-module:
\begin{equation}\label{eq:4.3}
  M_{A_n}(2j)=\bigoplus_{0\leq 2k\leq n} M_{A_{n-1}}(2k)\tensor L(c_n,h_{2k+1,2j+1}^{(n)}).
\end{equation}
Let $J$ be the maximal ideal of $V$.
Then $V/J\cong M_{A_n}$ by Theorem \ref{thm:3.6}.
Since $V^{[n-1]}\tensor \la \eta^n\ra$ is rational, there is a
$V^{[n-1]}\tensor \la \eta^n\ra$-submodule $X$ such that $V=X\oplus J$.
It follows from \eqref{eq:4.3} that $X$ is isomorphic to 
\[
  X \cong \bigoplus_{0\leq 2k\leq n}M_{A_{n-1}}(2k)\tensor L(c_n,h_{2k+1,1}^{(n)})
\]
as a $V^{[n-1]}\tensor \la \eta^n\ra$-module.
Since the Griess algebra of $V$ and that of $M_{A_n}$ are isomorphic,
we have $J\cap V_2=0$ and $V_2\subset X$.
Therefore, it suffices to show that $X$ forms a subalgebra of $V$ since $V$ is generated by
the Griess algebra.
By Lemma \ref{lem:4.4}, $J$ has a composition series, and since $V$ is of CFT-type,
there is no composition factor isomorphic to $M_{A_n}$ in $J$.
It follows from the fusion rules in \eqref{eq:3.2} that 
$L(c_n,h_{2k+1,1}^{(n)})$, $0\leq 2k\leq n$ are closed under the fusion product.
On the other hand, there is no irreducible $\la \eta^n\ra$-submodule isomorphic to
$L(c_n,h_{2k+1,1}^{(n)})$ in $M_{A_n}(2j)$ if $j>0$ by \eqref{eq:4.3} since 
$h_{r,s}^{(n)}=h_{r',s'}^{(n)}$ if and only if $(r,s)=(r',s')$ or $(r+r',s+s')=(n+2,n+3)$ 
by \eqref{eq:3.1}.
Thus, $\la X\ra \cap J=0$.
This completes the proof.
\qed

\appendix
\section{Appendix}

Here we prove Lemma \ref{lem:3.9}.
Consider the generating series
\begin{equation}\label{eq:a.1}
  \psi(z)=\dsum_{n\in \Z} \psi_{n+1/2}z^{-n-1},~~~~~
  [\psi_r,\psi_s]_+=\delta_{r+s,0},~~~r,s\in \Z+1/2.
\end{equation}
Then the SVOA $L(\shf,0)\oplus L(\shf,\shf)$ can be explicitly realized as follows (cf.~\cite{KR}).
\[
\begin{array}{l}
  L(\shf,h)
  =\Span_\C\{ \psi_{-r_1}\cds \psi_{-r_k} \vac \mid r_1>\cds >r_k> 0,~ k\equiv 2h ~\bmod 2\},~~
  h=0,1/2,
  \medskip\\
  Y(\vac,z)=\id,~~~~~~
  Y(\psi_{-1/2}\vac,z)=\psi(z),~~~~~~
  \w=\dfr{1}{2}\psi_{-3/2}\psi_{-1/2}\vac,
\end{array}
\]
where $\psi_r \vac=0$ for $r>0$.

Let $e\in V$ be an Ising vector. Let $b\in V_2$ be a highest weight vector for $\la e\ra$
with highest weight $1/2$ and $y\in V_N$  a highest weight vector for $\la e\ra$
with highest weight $h$.
Consider $V$ as a $\la e\ra\tensor \com_V \la e\ra$-module.
We can write $2e=\psi_{-3/2}\psi_{-1/2}\vac$, $b=\psi_{-1/2}\vac \tensor v$ and
$y=w\tensor x$, where $\psi_{-1/2}\vac$ and $w$ are highest weight vectors of
$L(\shf,\shf)$ and $L(\shf,h)$, respectively.
By \eqref{eq:3.2}, $L(\shf,0)$ and $L(\shf,\shf)$ are simple current modules. 
It then follows from Theorem 2.10 of \cite{ADL} that we have a factorization
\begin{equation}\label{eq:a.2}
  Y(b,z)y=\psi(z) w\tensor J(v,z)x
\end{equation}
where $J(\cd,z)$ is an intertwining operator
among $\com_V\la e\ra$-submodules of $V$.

First, we consider the case $h=0$.
Then $w=\vac$ and as a $\la e\ra\tensor \com_V\la e\ra$-module,
we have 
\[
  e=\dfr{1}{2}\psi_{-3/2}\psi_{-1/2}\vac,~~~
  b=\psi_{-1/2}\vac \tensor v= \psi_{-1/2}v_{(-1)}\vac,~~~
  y= \vac \tensor x,~~~
  Y(b,z)=\psi(z)\tensor J(v,z).
\]
For $n>0$, we also have 
\[
  \psi_{-n-1/2}\vac =\dfr{1}{n!}{e_{(0)}}^n\,\psi_{-1/2}\vac .
\]
Since $b_{(k)}y\in V_{N-k+1}[\shf]_e$ (cf.~Eq.~\eqref{eq:3.3}), we have $b_{(k)}y = 0$ and
$v_{(k)}x=0$ if $k\geq N$.
(Note that $V_{\leq 1}[\shf]_e=0$ since $V$ is of OZ-type.)
We set $P_{1/2,k,j}=0$ if $k\geq N$ or $k+j\geq N$.
Let $i\geq 0$.
Then from \eqref{eq:a.2} with $w=\vac$, we have
\[
\begin{array}{ll}
  b_{(N-1)}y
  &=\psi_{-1/2}v_{(N-1)} x,
  \\
  b_{(N-2)}y
  &= \psi_{-1/2}v_{(N-2)} x+e_{(0)}\psi_{-1/2}v_{(N-1)} x,
  \\
  b_{(N-3)}y
  &= \psi_{-1/2}v_{(N-3)}x+e_{(0)}\psi_{-1/2}v_{(N-2)}x
  +\dfr{{e_{(0)}}^2}{2!}\psi_{-1/2}v_{(N-1)}x,
  \\
  &\vdots
  \\
  b_{(N-i-1)} y
  &= \dsum_{j=0}^i \dfr{{e_{(0)}}^j}{j!}\psi_{-1/2}v_{(N-i+j-1)}x.
\end{array}
\]
Multiplying ${e_{(0)}}^{i-j}$ with $b_{(N-j-1)}y$ for $0\leq j\leq i$,
we obtain the following linear system:
\[
\begin{array}{ll}
  {e_{(0)}}^i\, b_{(N-1)}y
  &= {e_{(0)}}^i\, \psi_{-1/2}v_{(N-1)} x,
  \\
  {e_{(0)}}^{i-1}\, b_{(N-2)}y
  &= {e_{(0)}}^i \psi_{-1/2}v_{(N-1)} x +{e_{(0)}}^{i-1} \psi_{-1/2}v_{(N-2)} x,
  \\
  {e_{(0)}}^{i-2}\, b_{(N-3)}y
  &= \dfr{{e_{(0)}}^i}{2!}\psi_{-1/2}v_{(N-1)}x+ {e_{(0)}}^{i-1}\psi_{-1/2}v_{(N-2)}x
  +{e_{(0)}}^{i-2} \psi_{-1/2}v_{(N-3)}x,
  \\
  &\vdots
  \\
  b_{(N-i-1)} y
  &= \dsum_{j=0}^i \dfr{{e_{(0)}}^j}{j!}\psi_{-1/2}v_{(N-i+j-1)}x.
\end{array}
\]
This system can be formatted as follows.
\[
  \begin{bmatrix}
    {e_{(0)}}^i\,b_{(N-1)}y
    \medskip\\
    {e_{(0)}}^{i-1}\, b_{(N-2)}y
    \medskip\\
    {e_{(0)}}^{i-2}\, b_{(N-3)}y
    \medskip\\
    \vdots
    \medskip\\
    {e_{(0)}}b_{(N-i)}y
    \medskip\\
    b_{(N-i-1)}y
  \end{bmatrix}
  = \exp(\mathbb{J})
  \begin{bmatrix}
    {e_{(0)}}^i \psi_{-1/2}v_{(N-1)}x
    \medskip\\
    {e_{(0)}}^{i-1} \psi_{-1/2}v_{(N-2)}x
    \medskip\\
    {e_{(0)}}^{i-2} \psi_{-1/2}v_{(N-3)}x
    \medskip\\
    \vdots
    \medskip\\
    {e_{(0)}} \psi_{-1/2}v_{(N-i)}x
    \medskip\\
    \psi_{-1/2}v_{(N-i-1)}x
  \end{bmatrix}_,
  ~~~
  \mathbb{J}=
  \begin{bmatrix}
    ~0& & & & &
    \\
    ~1 & 0 & & & &
    \\
    & 1 & 0 & & &
    \\
    & &\ddots &\ddots &  &
    \\
    & & & 1 & 0 &
    \\
    & & & & 1 & 0~
  \end{bmatrix}_.
\]
Therefore, we obtain
\[
  \begin{bmatrix}
    {e_{(0)}}^i\, \psi_{-1/2}v_{(N-1)}x
    \medskip\\
    {e_{(0)}}^{i-1}\, \psi_{-1/2}v_{(N-2)}x
    \medskip\\
    {e_{(0)}}^{i-2}\, \psi_{-1/2}v_{(N-3)}x
    \medskip\\
    \vdots
    \medskip\\
    {e_{(0)}} \psi_{-1/2}v_{(N-i)}x
    \medskip\\
    \psi_{-1/2}v_{(N-i-1)}x
  \end{bmatrix}
  = \exp(-\mathbb{J})
  \begin{bmatrix}
    {e_{(0)}}^ib_{(N-1)}y
    \medskip\\
    {e_{(0)}}^{i-1} b_{(N-2)}y
    \medskip\\
    {e_{(0)}}^{i-2} b_{(N-3)}y
    \medskip\\
    \vdots
    \medskip\\
    {e_{(0)}}b_{(N-i)}y
    \medskip\\
    b_{(N-i-1)}y
  \end{bmatrix}_,
\]
and
\[
  \psi_{-1/2}v_{(N-i-1)}x
  = \dsum_{j=0}^i \dfr{(-e_{(0)})^{j}}{j!}b_{(N-i+j-1)}y
  = b_{(N-i-1)}y+\dsum_{j=1}^i \dfr{(-e_{(0)})^{j}}{j!}b_{(N-i+j-1)}y.
\]
The left hand side is a highest weight vector for $\la e\ra$ and we obtain the lemma by  setting 
$P_{1/2,N-i-1,j}=(-e_{(0)})^j/j!$.

Next we consider the case $h=1/2$ and $w=\psi_{-1/2}\vac$.
Since $\psi(z)\psi_{-1/2}\vac \in \la e\ra$, there exist $Q_{m+2}\in \mathrm{U}(\vir(e)_-)[m+2]$
such that $\psi_{-m-3/2}\psi_{-1/2}\vac = Q_{m+2} \vac$
and
\begin{equation}\label{eq:a.3}
  \psi(z)\psi_{-1/2}\vac
  = \vac z^{-1}+\dsum_{m\geq 0}\psi_{-m-3/2}\psi_{-1/2}\vac z^{m+1}
  = \vac z^{-1}+\dsum_{m\geq 0} Q_{m+2} \vac z^{m+1}.
\end{equation}
Note that $Q_{m+2}\in \mathrm{U}(\vir(e)_-)[m+2]$ is not unique but $Q_{m+2}\vac$ is uniquely
determined in $\la e\ra\cong L(\shf,0)$.
Since $b_{(k)}y\in V_{N-k+1}=0$ if $k>N+1$, $b_{(N+1)}y\in V_0=\C \vac$ and $b_{(N)}y\in V_1=0$,
we set $P_{0,k,j}=0$ for $k\geq N$ and $j\geq 0$.
We also set $P_{0,k,j}=0$ if $k+j>N+1$.
We will define $P_{0,k,j}$ with $j\geq 0$ recursively for $k<N$ as follows.

Suppose $i>0$ and for $k> N-i$ there exist $P_{0,k,j}$, $j\geq 0$, such that
\begin{equation}\label{eq:a.4}
  v_{(k-1)}x=b_{(k)}y+\dsum_{j>0} P_{0,k,j} b_{(k+j)}y.
\end{equation}
By \eqref{eq:a.2} and \eqref{eq:a.3} we have
\[
  b_{(N-i)}y=v_{(N-i-1)}x +\dsum_{0\leq m\leq i-1} Q_{m+2} \,v_{(N-i+m+1)} x
\]
and
\[
  v_{(N-i-1)}x=b_{(N-i)}y-\dsum_{0\leq m\leq i-1} Q_{m+2}
   \l( b_{(N-i+m+2)}y +\dsum_{j>0}P_{0,N-i+m+2,j} b_{(N-i+m+2+j)}y\r) .
\]
Expanding the right hand side, we obtain recursive relations for
$P_{0,N-i,j}$ such that \eqref{eq:a.4} holds for $k=N-i$.
Note that all $P_{0,k,j}$ are determined only by $k$, $j$, $b_{(N+1)}y\in V_0$ and
$Q_{m+2}\vac \in \la e\ra$, and $Q_{m+2}\vac$ is uniquely determined only by the
structure of $L(\shf,0)$ and independent of the structure of $V$ itself.
This completes the proof of Lemma \ref{lem:3.9}.
\qed

\small

\end{document}